\documentclass[12pt,reqno]{amsart}

\usepackage{amsmath}
\usepackage{amssymb}
\usepackage{amscd}
\usepackage{amsfonts}

\newtheorem{theorem}{Theorem}[section]
\newtheorem{lemma}[theorem]{Lemma}
\newtheorem{proposition}[theorem]{Proposition}
\newtheorem{definition}[theorem]{Definition}
\newtheorem{remark}[theorem]{Remark}
\newtheorem{corollary}[theorem]{Corollary}
\newcommand{\vol}{\mathrm{Vol}}

\numberwithin{equation}{section}

\begin{document}

\title{Davies Type Estimate And The Heat Kernel Bound Under The Ricci Flow}
\author{Meng Zhu}
\thanks{Research is partially supported by China Postdoctral Science Foundation Grant No. 2013M531105}
\address{Shanghai Center for Mathematical Sciences\\Fudan University\\220 Handan Road, Shanghai 200433, China}
\email{mengzhu@fudan.edu.cn}
\date{}

\begin{abstract}
We prove a Davies type double integral estimate for the heat kernel $H(y,t;x,l)$ under the Ricci flow. As a result, we give an affirmative answer to a question proposed in \cite{Chowetc2010}. Moreover, we apply the Davies type estimate to provide a new proof of the Gaussian upper and lower bounds of $H(y,t;x,l)$ which were first shown in \cite{CTY2011}.
\end{abstract}

\maketitle

\section{Introduction}

On a complete Riemannian manifold $(M^n, g_{ij})$, the heat kernel $H(x,y,t)$, is the smallest positive fundamental solution to the heat equation
\begin{equation}\label{he}
\frac{\partial u}{\partial t}=\Delta u.
\end{equation}
Heat kernel estimates is of great importance and interest due to its relation to many other properties of the manifolds, such as Harnack estimate, Sobolev inequality, Log Soblev inequality, Faber-Krahn type inequality, and Nash type inequality (see e.g. \cite{Var1985}, \cite{Dav1987}, \cite{Gri1994}, \cite{CKS1987}, \cite{LiYa1986}, \cite{Car1996}). Since the work of Nash \cite{Nas1958} and Aronson \cite{Aro1967}, many methods have been discovered for deriving Gaussian upper and lower bounds of $H(x,y,t)$, see e.g., \cite{CLY1981}, \cite{LiYa1986}, \cite{Dav1989}, \cite{Gri1997}, \cite{Dav1992}, \cite{LiWa1999}. One of the methods was developed by Li-Wang in \cite{LiWa1999}. They obtained a Gaussian upper bound for $H(x,y,t)$ based on the parabolic mean value inequality and the following double integral upper estimate of the heat kernel proved by E.B. Davies \cite{Dav1992}:

\begin{theorem} Let $(M, g)$ be a complete Riemannian manifold. For any two bounded subsets $U_1$ and $U_2$ of $M$, one has:
\begin{equation}
\int_{U_1}\int_{U_2}H(x,y,t)d\mu(x)d\mu(y)\leq \vol^{\frac{1}{2}}(U_1)\vol^{\frac{1}{2}}(U_2)e^{-\frac{d^2(U_1, U_2)}{4t}},
\end{equation}
where $d(U_1,U_2)$ is the distance between $U_1$ and $U_2$.
\end{theorem}

In this paper, we consider the heat kernel of the time-evolving heat equation under the Ricci flow on a complete manifold $M^n$, i.e.,
\begin{equation}\label{tehe}
\frac{\partial u}{\partial t}=\Delta_t u,
\end{equation}
where $\Delta_t$ is the Laplacian with respect to a complete solution $g_{ij}(t)$, $t\in [0,T)$ and $T<\infty$, of the following Ricci flow equation
\begin{equation}\label{rf}
\left\{\begin{aligned}
&\frac{\partial g_{ij}(t)}{\partial t}=-2R_{ij}\\
& g(0)=g_0
\end{aligned}
\right.
\end{equation}
on $M$.

The existence and uniqueness of the heat kernel $H(y,t;x,l)$ to \eqref{tehe} were proved in \cite{Gue2002} and \cite{Chowetc2010}. When $M$ is compact, C. Guenther \cite{Gue2002} obtained a Gaussian lower bound of $H(y,t;x,l)$. For general complete manifolds, L. Ni \cite{Ni2004} first obtained Gaussian estimates of $H(y,t;x,l)$ assuming uniformly bounded curvature and nonnegative Ricci curvature along the Ricci flow. Q. Zhang \cite{Qzhang2010} proved both Gaussian upper and lower bounds of $H(y,t;x,l)$ for type I ancient $\kappa$-solutions of the Ricci flow with nonnegative curvature operator. In \cite{Xu2012}, G. Xu removed the assumption on the nonnegativity of the curvature operator, and obtained similar Gaussian bounds of the heat kernel for general type I ancient $\kappa$-solutions of the Ricci flow. Assuming that $M$ has nonnegative Ricci curvature and is not Ricci flat for all $t$, Gaussian estimates of $H(y,t;x,l)$ were given by Cao-Zhang \cite{CaZh2011}. When the Ricci curvature is only assumed to be uniformly bounded, Chau-Tam-Yu \cite{CTY2011} also showed Gaussian lower and upper bounds for $H(y,t;x,l)$, but the constants in the estimates depend on more geometric quantities than the ones in Cao-Zhang's estimates (see \cite{Chowetc2010} for more detail). While the authors in \cite{CTY2011} and \cite{Chowetc2010} mainly used Grigor'yan's method in \cite{Gri1997} to prove the Gaussian upper bound of the heat kernel, they raised the following question:\\

\hspace{-.45cm}\textbf{Question:} Is there a "Ricci flow" version of Davies' estimate as the one in Theorem 1.1?\\

Our main goal is to give an affirmative answer to the question above. More specifically, we prove

\begin{theorem}\label{davies}
Let $(M^n, g_{ij}(t))$ be a complete solution to \eqref{rf} on $[0,T)$ and $T<\infty$. Suppose that $H(y,t;x,l)$ is the heat kernel of \eqref{tehe}, and $Rc\geq -K_1$ on $[0,T)$ for some nonnegative constant $K_1$. Then for any open sets, $U_1$ and $U_2$, with compact closure in $M$, and $0\leq l<t<T$, we have
\begin{align*}
\int_{U_1}\int_{U_2}H(y,t;x,l)d\mu_t(y)d\mu_l(x)&\leq e^{C_0K_1T}e^{-\frac{d^2_t(U_1,U_2)}{4e^{4K_1T}(t-l)}}\vol_l(U_1)^{\frac{1}{2}}\vol_t(U_2)^{\frac{1}{2}},
\end{align*}
where $C_0$ is a constant depends only on $n$.
\end{theorem}

Using Theorem \ref{davies}, we are able to provide a new proof of the Gaussian upper bound and lower bounds of $H(y,t;x,l)$ which was first shown by Chau-Tam-Yu \cite{CTY2011}. For the Gaussian upper bound, we have

\begin{theorem}\label{upperbound}
Under the Ricci flow, assume that $Rc\geq-K_1$ on $M\times [0, T)$ for some nonnegative constant $K_1$ and $T<\infty$, and that $\Lambda=\int_0^T\sup_M|Rc|(t)dt<\infty$. We have the following upper bound:
\begin{align*}
H(y,t;x,l)\leq C_1e^{C_2e^{C_3\Lambda+C_4K_1T}} \min\left\{\frac{exp\left(-\frac{d^2_{t}(x,y)}{8e^{4K_1T}(t-l)}\right)}{\vol_{l}(B_{l}(y,\sqrt{\frac{t-l}{8}}))},\frac{exp\left(-\frac{d^2_{t}(x,y)}{8e^{4K_1T}(t-l)}\right)}{\vol_{t}(B_{t}(x,\sqrt{\frac{t-l}{8}}))}\right\},
\end{align*}
where constants $C_1,\ C_2,$, $C_3$ and $C_4$ depend only on $n$.
\end{theorem}

Next, following a method of Li-Tam-Wang \cite{LTW1997}, we can use a gradient estimate of Q. Zhang in \cite{Zha2006} (see also \cite{CaHa2009}) to show a Gaussian lower bound of $H(y,t;x,l)$.

\begin{theorem}\label{lowerbound}
Let $(M^n, g_{ij}(t))$, $t\in[0, T)$ and $T<\infty$, be a complete solution to the Ricci flow \eqref{rf}. Suppose that $Rc\geq-K_1$ on $M\times[0,T)$ and $\Lambda=\int_0^T\sup_M|Rc|(t)dt<\infty$. Then we have the following lower bound:
\begin{align*}
& H(y,t;x,l)\geq\frac{C_5e^{-C_6e^{C_7\Lambda+C_8K_1T}} exp\left(\displaystyle -\frac{4d^2_{t}(x,y)}{(t-l)}\right)}{\vol_{t}(B_{t}(x, \sqrt{\frac{t-l}{8}}))},
\end{align*}
where $C_5,\ C_6$, $C_7$ and $C_8$ are positive constants only depending on $n$.
\end{theorem}

This paper is organized as follows. In section 2, we review some known facts about the fundamental solutions of heat-type equations with metric evolving under a group of more general equations than the Ricci flow, including the existence, uniqueness and the mean value inequality. In section 3, we prove Theorem \ref{davies}. The main idea is similar to that in the fixed metric case. However, since the heat kernel is not self-symmetric in this case, we use the symmetry between the heat kernel and the adjoint heat kernel instead. In section 4, we first use the method in \cite{LiWa1999} to prove a slightly different version of the Gaussian upper estimate from the one in Theorem \ref{upperbound}. Then we show that Theorem \ref{upperbound} can be derived from this upper estimate and an $L^1$ bound of $H$. In section 5, we finish the proof of Theorem \ref{lowerbound}.\\

\hspace{-.45cm}\textbf{Acknowledgements:} The author would like to thank Professor Huai-Dong Cao for many valuable suggestions, and for his constant support and encouragement. I also want to thank Professors Qing Ding, Jixiang Fu, Jiaxing Hong, Hong-Quan Li, Jun Li, Jiaping Wang, Quanshui Wu, Guoyi Xu and Weiping Zhang for helpful conversations and encouragement, and Professors Luen-Fai Tam and Qi S. Zhang for their interests in this work. In addition, I want to express my appreciation to the Shanghai Center for Mathematical Sciences, where this work was carried out, for its hospitality and support.

\vspace{.2cm}

\section{Preliminaries}

In this section, we present some basic results regarding the heat kernel and adjoint heat kernel of time evolving heat-type equations. The readers can refer to \cite{Chowetc2010} for more detail.

Let $M^n$ be a Riemannian manifold, and $g_{ij}(t)$, $t\in [0, T)$ and $T<\infty$, a complete solution to the following equation:
\begin{equation} \label{flow}
\left\{\begin{aligned}
&\frac{\partial g_{ij}(t)}{\partial t}=-2A_{ij}\\
& g(0)=g_0,
\end{aligned}
\right.
\end{equation}
where $A_{ij}(t)$ is a time-dependent symmetric 2-tensor. Consider the following heat operator with potential:
$$L=\frac{\partial }{\partial t}-\Delta_t + Q,$$
where $\Delta_t$ is the Laplacian with respect to $g_{ij}(t)$ and $Q:\ M\times[0, T)\rightarrow \mathbb{R}$ is a $C^{\infty}$ function.

\begin{definition} Let $\mathbb{R}_T^2=\{(t,l)\in\mathbb{R}^2\ |\ 0\leq l<t<T\}$. A fundamental solution for the operator $L$ is a function
$$H: M\times M\times \mathbb{R}_T^2\rightarrow \mathbb{R}$$
that satisfies:
\begin{itemize}
\item[(1)] $H$ is continuous, $C^2$ in the first two space variables and $C^1$ in the last two time variables,
\item[(2)] $LH=(\frac{\partial }{\partial t}-\Delta_{t,y} + Q)H(y,t;x,l)=0$,
\item[(3)] $\lim_{t\searrow l}H(y,t;x,l)=\delta_x$,
\end{itemize}
where $H(y,t;x,l)=H(y,x,(t,l))$.

The heat kernel for $L$ is defined to be the minimal positive fundamental solution.
\end{definition}

\begin{definition}
The adjoint heat kernel is the minimal positive fundamental solution
$$G: M\times M\times \mathbb{R}_T^2\rightarrow \mathbb{R}$$
for the operator $L^*=\frac{\partial }{\partial l} +\Delta_l-Q-A$, i.e., $G$ satisfies
\begin{itemize}
\item[(1)] $G$ is continuous, $C^2$ in the first two space variables and $C^1$ in the last two time variables,
\item[(2)] $L^*G=(\frac{\partial }{\partial l}+\Delta_{l,x} - Q-A)G(x,l;y,t)=0$,
\item[(3)] $\lim_{l\nearrow t}G(x,l;y,t)=\delta_y$,
\end{itemize}
where $G(x,l;y,t)=G(y,x,(t,l))$, and $A=tr_{g(l)}(A_{ij})$.
\end{definition}

First of all, we have the following existence and uniqueness of the heat kernel and the adjoint heat kernel according to \cite{Gue2002} and \cite{Chowetc2010}.
\begin{theorem}
Let $M^n$ be a complete manifold, and $g(t)$, $t\in[0, T)$ for some $T<\infty$, a smooth family of Riemannian metrics on $M$. If $ \int_0^T -\displaystyle \inf_MQ(t)dt<\infty$, then there exists a unique $C^{\infty}$ minimal positive fundamental solution $H(y,t;x,l)$ for the operator $L=\frac{\partial }{\partial t}-\Delta_{t,y} + Q$.
\end{theorem}

As in the fixed metric case, the heat kernel $H(y,t;x,l)$  for $L$ on $M$ is the limit of the Dirichlet heat kernels on a sequence of exhausting subsets in $M$.

\begin{definition}
Let $\Omega\subset M $ be a bounded subset, the Dirichlet Heat Kernel on $\Omega$ for $L$, denoted by $H_{\Omega}(y,t;x,l)$ is the fundamental solution to $\frac{\partial u}{\partial t}=\Delta_{t} u- Q u$ in $\Omega$ and satisfies

i) $\displaystyle \lim_{t\searrow l}H_{\Omega}=\delta_x$ for $x \in \textrm{Int} (\Omega)$;

ii) $H_{\Omega}(y,t;x,l)=0$, for $y\in \partial\Omega$ and $x\in \textrm{Int}(\Omega)$.
\end{definition}

\begin{proposition}[see e.g. \cite{Chowetc2010}]\label{prop1}
Let $\Omega_i\subset M $ be a sequence of exhausting bounded sets, and $H_{\Omega_i}(y,t;z,s)$ the Dirichlet Heat Kernel on $\Omega_i$. Then
$$\lim_{i\rightarrow \infty}H_{\Omega_i}(y,t;z,s)=H(y,t;z,s)$$
uniformly on any compact subset of $M\times M\times \mathbb{R}_T^2$.
\end{proposition}

Moreover, the heat kernel and the adjoint heat kernel satisfy the following important properties:

\begin{proposition}[see e.g. \cite{Chowetc2010}]\label{prop2}
With the assumptions above, we have
\begin{itemize}
\item[(1)] $H_{\Omega}(y,t;x,l)=G_{\Omega}(x,l;y,t)$, for $x, y\in \Omega$;
\item[(2)] $H(y,t;x,l)=G(x,l;y,t)$;
\item[(3)] $H_{\Omega}(y,t;x,l)=\int_{\Omega} H_{\Omega}(y,t;z,s)H_{\Omega}(z,s;x,l)d\mu_{g(s)}(z)$;
\item[(4)] $H(y,t;x,l)=\int_M H(y,t;z,s)H(z,s;x,l)d\mu_{g(s)}(z)$;
\end{itemize}
where $\Omega$ is an open subset in $M$ with compact closure, and $G_{\Omega}(x,l;y,t)$ is the Dirichlet heat kernel for the adjoint operator $\frac{\partial }{\partial l} +\Delta_l-Q-A$.
\end{proposition}

Let $$\Lambda=\int_0^T\sup_{M}|A_{ij}|_{g(t)}(t)dt.$$ Assume that there exists a metric $g^{\prime}$ and a positive constant $\hat{C}\geq 1$ such that
$$\hat{C}^{-1}g^{\prime}\leq g(0)\leq \hat{C}g^{\prime},$$
then we have
$$\hat{C}^{-1}e^{-2\Lambda}g^{\prime}\leq g(t)\leq \hat{C}e^{2\Lambda}g^{\prime}$$
for all time $t\in[0, T)$.

Suppose that $u:M\times[0,T)\rightarrow \mathbb{R}$ is a positive subsolution to
$$\frac{\partial u}{\partial t}\leq \Delta_{t}u-Qu.$$
Define the parabolic cylinder:
$$P_{g^{\prime}}(x,\tau,r,-r^2)=B_{g^{\prime}}(x,r)\times[\tau-r^2,\tau],$$
where $B_{g^{\prime}}(x,r)$ represents the geodesic ball of radius $r$ centered at $x$ in $M$ with respect to the metric $g^{\prime}$.

By Moser iteration, Chau-Tam-Yu \cite{CTY2011} got the following mean value inequality (see also \cite{Chowetc2010}):
\begin{theorem}\label{mean} In the above setting, assume that $Rc(g^{\prime})\geq-K_1$ on $M$ with some $K_1\geq 0$. Then
$$\sup_{P_{g^{\prime}}(x_0,t_0,r_0,-(r_0)^2)}u\leq \frac{C_1\hat{C}^{\frac{n(n+3)}{2}}e^{C_2\Lambda+C_3\sqrt{K_1}r_0+ \check{C}t_0}}{r_0^2\vol_{g^{\prime}}(B_{g^{\prime}}(x_0,r_0))}\int\limits_{P_{g^{\prime}}(x_0,t_0,2r_0,-(2r_0)^2)}u(x,s)d\mu_{g^{\prime}}(x)ds,$$
where $C_1$, $C_2$ and $C_3$ are constants only depending on $n$, and $\check{C}=-\inf_{M\times[0,T)}\{Q+\frac{1}{2}A\}$. 
\end{theorem}

\section{Davies type estimate for the heat kernel under the Ricci flow}

Assume that $(M^n, g_{ij}(t))$ is a complete solution to the Ricci flow \eqref{rf}
for $t\in[0,T)$ and $T<\infty$. 

Denote by $H(y,t;x,l)>0$, $0\leq l<t<T$, the Heat Kernel of \eqref{tehe} under the Ricci flow, i.e.,
\begin{equation}\label{hk}
\left\{ \begin{aligned}
&\frac{\partial H}{\partial t}=\Delta_{t,y} H\\
&\displaystyle \lim_{t\searrow l}H=\delta_x.
\end{aligned}\right.
\end{equation}

Then $G(x,l;y,t)=H(y,t;x,l)$ is the adjoint Heat Kernel to the following conjugate heat equation:
\begin{equation}\label{ahk}
\left\{ \begin{aligned}
&\frac{\partial G}{\partial l}=-\Delta_{l,x} G + RG\\
&\displaystyle \lim_{l\nearrow t}G=\delta_y,
\end{aligned}\right.
\end{equation}
where $R$ is the scalar curvature of $M$.

In the remaining of this paper, we will use the following notations: $B_t(x,r):=B_{g(t)}(x,r)$, $\vol_s(U):=\vol_{g(s)}(U)$, $d\mu_s:=d\mu_{g(s)}$, where $U$ is a subset of $M$ and $d\mu_{g(s)}$ denotes the volume element of $g(s)$.\\


\begin{proof}[Proof of Theorem \ref{davies}]: Let $\Omega_i\subset M $ be a sequence of exhausting bounded sets, and $H_{\Omega_i}(y,t;z,s)$ the Dirichlet Heat Kernel on $\Omega_i$. Since $U_1$ and $U_2$ are bounded, we may assume that $(U_1\bigcup U_2)\subset\Omega_i$ for each $i$.

By Propositions \ref{prop1} and \ref{prop2}, we have
\begin{align*}
&\int_{U_1}\int_{U_2}H(y,t;x,l)d\mu_t(y)d\mu_l(x)\\
=&\displaystyle \lim_{i\rightarrow\infty}\int_{U_1}\int_{U_2}\int_{\Omega_i}H_{\Omega_i}(y,t;z,\frac{t+l}{2})H_{\Omega_i}(z,\frac{t+l}{2};x,l)d\mu_{\frac{t+l}{2}}(z)d\mu_t(y)d\mu_l(x)\\
=&\displaystyle \lim_{i\rightarrow\infty}\int_{\Omega_i}\int_{U_1}\int_{U_2}H_{\Omega_i}(y,t;z,\frac{t+l}{2})H_{\Omega_i}(z,\frac{t+l}{2};x,l)d\mu_t(y)d\mu_l(x)d\mu_{\frac{t+l}{2}}(z)\\
=&\lim_{i\rightarrow\infty}\int_{\Omega_i}u_i(z,\frac{t+l}{2})v_i(z,\frac{t+l}{2})d\mu_{\frac{t+l}{2}}(z),
\end{align*}
where $$u_i(z,s)=\int_{U_1}H_{\Omega_i}(z,s;x,l)d\mu_l(x),$$
and $$v_i(z,s)=\int_{U_2}H_{\Omega_i}(y,t;z,s)d\mu_t(y)=\int_{U_2}G_{\Omega_i}(z,s;y,t)d\mu_t(y).$$ Here $G_{\Omega_i}(z,s;y,t)$ denotes the adjoint Dirichlet Heat Kernel on $\Omega_i$.

From $Rc\geq -K_1$, we have for any $s\in[0, T)$,
$$\ \ g_{ij}(s)\leq e^{2K_1T}g_{ij}(0),$$
$$\textrm{and}\quad d_s(x,y)\leq e^{K_1T}d_0(x,y),$$
where $d_s(x,y)$ is the distance function at time $s$.

Let $\xi(z,s)=\frac{d^2_0(z,U_1)}{2C_{K_1T}(s-l)}$ for $C_{K_1T}=e^{2K_1T}$ and $s>l$. Since $$|\nabla d_0(z,U_1)|_{g(s)}\leq e^{K_1T}|\nabla d_0(z,U_1)|_{g(0)}\leq e^{K_1T},$$ we have
$$\frac{\partial \xi}{\partial s}+\frac{1}{2}|\nabla\xi|^2_{g(s)}\leq -\frac{d^2_0(z,U_1)}{2C_{K_1T}(s-l)^2}+\frac{d^2_0(z,U_1)\cdot e^{2K_1T}}{2C_{K_1T}^2(s-l)^2}=0.$$

We compute
\begin{align*}
&\frac{d}{ds}\int_{\Omega_i}u^2_i(z,s)e^{\xi(z,s)}d\mu_s(z)\\
=&\int_{\Omega_i}(2u_i\Delta_{s,z}u_i-Ru^2_i+u^2_i\cdot\frac{\partial \xi(z,s)}{\partial s})e^{\xi(z,s)}d\mu_s(z).
\end{align*}

Notice that
\begin{align*}
2\int_{\Omega_i}u_i\Delta_{s,z}u_ie^{\xi}d\mu_s(z)&=-2\int_{\Omega_i}(|\nabla u_i|^2+u_i\nabla u_i\cdot\nabla\xi)e^{\xi}d\mu_s(z)+2\int_{\partial\Omega_i}u_i\partial_{\nu}u_ie^{\xi}dS\\
&\leq\int_{\Omega_i}\frac{1}{2}u^2_i|\nabla \xi|^2e^{\xi}d\mu_s(z),
\end{align*}
where we have used the Cauchy-Schwartz inequality for $-u_i\nabla u_i\cdot\nabla\xi$ and the fact that $u_i|_{\partial\Omega_i}=0$.

Thus,
$$\frac{d}{ds}\int_{\Omega_i}u^2_i(z,s)e^{\xi(z,s)}d\mu_s(z)\leq C_1K_1\int_{\Omega_i}u^2_i(z,s)e^{\xi(z,s)}d\mu_s(z),$$
where $C_1$ is a constant only depending on $n$.
Hence, we have
\begin{align*}
\int_{\Omega_i}u^2_i(z,s)e^{\xi(z,s)}d\mu_s(z)&\leq e^{C_1K_1(s-l)}\lim_{h\rightarrow l}\int_{\Omega_i}u^2_i(z,h)e^{\xi(z,h)}d\mu_h(z)\leq e^{C_1K_1(t-l)}\vol_l(U_1).
\end{align*}

Let $\eta(z,s)=\frac{d^2_0(z,U_2)}{2C_{K_1T}(t-s)}$ for $C_{K_1T}=e^{2K_1T}$, then
$$\frac{\partial\eta}{\partial s}-\frac{1}{2}|\nabla\eta|^2_{g(s)}\geq\frac{d^2_0(z, U_2)}{2C_{K_1T}(t-s)^2}-\frac{d^2_0(z,U_2)\cdot e^{2K_1T}}{2C_{K_1T}^2(t-s)^2}= 0.$$

Moreover,
\begin{align*}
&\frac{d}{ds}\int_{\Omega_i}v^2_i(z,s)e^{\eta(z,s)}d\mu_s(z)\\
=&\int_{\Omega_i}(-2u_i\Delta_{s,z}v_i+v^2_i\cdot\frac{\partial \eta(z,s)}{\partial s})e^{\eta(z,s)}d\mu_s(z).
\end{align*}

Since,
\begin{align*}
-2\int_{\Omega_i}v_i\Delta^s_zv_ie^{\eta}d\mu_s(z)&=2\int_{\Omega_i}(|\nabla v_i|^2+v_i\nabla v_i\cdot\nabla \eta)e^{\eta}d\mu_s(z)-2\int_{\partial\Omega}v_i\partial_{\nu}v_ie^{\eta}dS\\
&\geq-\frac{1}{2}\int_{\Omega_i}v^2_i|\nabla\eta|e^{\eta}d\mu_s(z),
\end{align*}
we have
$$\frac{d}{ds}\int_{\Omega_i}v^2_i(z,s)e^{\eta(z,s)}d\mu_s(z)\geq 0.$$

It implies that, for large $i$,
\begin{align*}
\int_{\Omega_i}v^2_i(z,s)e^{\eta(z,s)}d\mu_s(z)&\leq \lim_{h\rightarrow t}\int_{\Omega_i}v^2_i(z,h)e^{\eta(z,h)}d\mu_h(z)\\
&=\vol_t(U_2).
\end{align*}

Now since
$$\frac{1}{2}\xi(z,\frac{t+l}{2})+\frac{1}{2}\eta(z,\frac{t+l}{2})\geq\frac{d^2_0(U_1,U_2)}{4C_{K_1T}(t-l)},$$
we have
\begin{align*}
&e^{\frac{d^2_0(U_1,U_2)}{4C_{K_1T}(t-l)}}\int_{U_1}\int_{U_2}H(y,t;x,l)d\mu_t(y)d\mu_l(x)\\
&=\lim_{i\rightarrow\infty}e^{\frac{d^2_0(U_1,U_2)}{4C_{K_1T}(t-l)}}\int_{\Omega_i}u_i(z,\frac{t+l}{2})v_i(z,\frac{t+l}{2})d\mu_{\frac{t+l}{2}}(z)\\
&\leq \lim_{i\rightarrow\infty}\int_{\Omega_i}u_i(z,\frac{t+l}{2})e^{\frac{1}{2}\xi(z,\frac{t+l}{2})}\cdot v_i(z,\frac{t+l}{2})e^{\frac{1}{2}\eta(z,\frac{t+l}{2})}d\mu_{\frac{t+l}{2}}(z)\\
&\leq e^{\frac{C_1K_1(t-l)}{2}}\vol_l(U_1)^{\frac{1}{2}}\vol_t(U_2)^{\frac{1}{2}},
\end{align*}

i.e.,
\begin{align*}
\int_{U_1}\int_{U_2}H(y,t;x,l)d\mu_t(y)d\mu_l(x)&\leq e^{-\frac{d^2_0(U_1,U_2)}{4C_{K_1T}(t-l)}}e^{\frac{C_1K_1(t-l)}{2}}\vol_l(U_1)^{\frac{1}{2}}\vol_t(U_2)^{\frac{1}{2}}.
\end{align*}
The Theorem follows from the fact that $d_t(U_1,U_2)\leq e^{K_1T}d_0(U_1,U_2)$.
\end{proof}

\vspace{.25cm}

\section{A Gaussian upper bound of $H(y,t;x,l)$}

Recall that by the Volume Comparison Theorem, we have the following Lemma (see e.g. \cite{Heb1999}):

\begin{lemma}\label{volumecomparison}
Let $(M^n, g_{ij})$ be a complete Riemannian manifold. If $Rc\geq -K_1$, for some constant $K_1\geq 0$, then for any point $x\in M$ and any $0<r\leq R$, we have
$$\vol(B(x, R))\leq (\frac{R}{r})^ne^{\sqrt{(n-1)K_1}R}\vol(B(x, r)),$$
where $B(x, R)$ is the geodesic ball of radius $R$ in $M$ centered at $x$.

In particular, letting $r\rightarrow 0$, we have
$$Vol(B(x, R))\leq CR^ne^{\sqrt{(n-1)K_1}R}$$
for some constant $C$ only depending on $n$.
\end{lemma}

Since we have obtained a Davies type estimate in Theorem \ref{davies}, similarly to the method in \cite{LiWa1999}, we can show a Gaussian upper bound of the heat kernel (see \cite{CTY2011}) by using the mean value inequality in Theorem \ref{mean}. In the following, unless otherwise stated, $C$, $C_0$, $C_1$, $C_2,\ \cdots$, and $\tilde{C}$, $\tilde{C}_0$, $\tilde{C}_1$, $\tilde{C}_2,\ \cdots$ all represent positive constants only depending on $n$.


\begin{theorem}\label{upperbound1}
Under the Ricci flow, assume that $Rc\geq -K_1$ on $M\times [0, T)$ and $\Lambda=\int_0^T\sup_M|Rc|(t)dt<\infty$. We have the following upper bound
\begin{align*}
& \ H(y,t;x,l)\leq \frac{C_{1}e^{C_2\Lambda+C_3 K_1T+C_4 \sqrt{K_1T}}exp\left(\displaystyle -\frac{d_{t}^2(x,y)}{8e^{4K_1T}(t-l)}\right)}{\sqrt{Vol_{l}(B_{l}(y,\sqrt{\frac{t-l}{8}}))}\sqrt{Vol_{t}(B_{t}(x,\sqrt{\frac{t-l}{8}}))}}.
\end{align*}
\end{theorem}

\begin{proof}
We have by Theorem \ref{mean} that, for $0\leq l_0<t_0< T$ and $r_0$ with $t_0-l_0\geq 4r_0^2$,
\begin{equation}\label{eq1}
\begin{aligned}
& \ H(y_0,t_0;x_0,l_0)\\
\leq &\ \frac{C_0e^{C_1\Lambda+C_2K_1T+C_3\sqrt{K_1}r_0}}{r_0^2Vol_{l_0}(B_{l_0}(y_0,r_0))}\int^{t_0}_{t_0-4r_0^2}\int_{B_{l_0}(y_0,2r_0)}H(y,s;x_0,l_0)d\mu_{l_0}(y)ds.\\
\end{aligned}
\end{equation}

Let $\tilde{l}=t_0-l$, $\tilde{s}=t_0-s$, and
$$\widetilde{G}(x, \tilde{l};y, \tilde{s})=G(x,t_0-\tilde{l};y,t_0-\tilde{s})=G(x,l;y,s)=H(y,s;x,l).$$

Then, $\widetilde{G}(x, \tilde{l};y, \tilde{s})$ satisfies
$$\frac{\partial \widetilde{G}(x,\tilde{l};y,\tilde{t})}{\partial \tilde{l}}=\widetilde{\Delta}_{\tilde{l},x} \widetilde{G} - \tilde{R}\,\widetilde{G}(x,\tilde{l};y,\tilde{t}),$$
where $\widetilde{\Delta}_{\tilde{l}}$ is the laplacian with respect to the solution $\tilde{g}(\tilde{l}):=g(t_0-\tilde{l})$ of the backward Ricci flow:
$$\frac{\partial \tilde{g}(\tilde{l})}{\partial \tilde{l}}=2\widetilde{R}_{ij}(\tilde{l}).$$
Thus, according to Theorem \ref{mean}, we have the following mean value inequality:
\begin{equation*}
\begin{aligned}
& \ \widetilde{G}(x_0, \tilde{l}_0;y, \tilde{s})\\
\leq &\ \frac{\tilde{C}_0e^{\tilde{C}_1\Lambda+\tilde{C}_2K_1T+\tilde{C}_3\sqrt{K_1}r_1}}{r_1^2Vol_{\tilde{g}(0)}(B_{\tilde{g}(0)}(x_0,r_1))}\int^{\tilde{l}_0}_{\tilde{l}_0-4r_1^2}\int_{B_{\tilde{g}(0)}(x_0,2r_1)}\widetilde{G}(x, \tilde{l};y, \tilde{s})d\mu_{\tilde{g}(0)}(x)d\tilde{l},
\end{aligned}
\end{equation*}
where $4r_1^2\leq t_0-l_0-4r_0^2\leq s-l_0$.

Hence Formula \eqref{eq1} becomes
\begin{equation*}
\begin{aligned}
& \ H(y_0,t_0;x_0,l_0)\\
\leq &\
\frac{C_0e^{C_1\Lambda+C_2K_1T+C_3\sqrt{K_1}r_0}}{r_0^2Vol_{l_0}(B_{l_0}(y_0,r_0))}\int^{4r_0^2}_0\int_{B_{l_0}(y_0,2r_0)}\widetilde{G}(x_0,\tilde{l}_0;y,\tilde{s})d\mu_{l_0}(y)d\tilde{s}\\
\leq &\
\frac{C_0e^{C_1\Lambda+C_2K_1T+C_3\sqrt{K_1}r_0}}{r_0^2Vol_{l_0}(B_{l_0}(y_0,r_0))}\frac{\tilde{C}_0e^{\tilde{C}_1\Lambda+\tilde{C}_2K_1T+\tilde{C}_3\sqrt{K_1}r_1}}{r_1^2Vol_{\tilde{g}(0)}(B_{\tilde{g}(0)}(x_0,r_1))}\\
\quad &\cdot\int^{4r_0^2}_0\int_{B_{l_0}(y_0,2r_0)}\int^{\tilde{l}_0}_{\tilde{l}_0-4r_1^2}\int_{B_{\tilde{g}(0)}(x_0,2r_1)}\widetilde{G}(x,\tilde{l};y, \tilde{s})d\mu_{\tilde{g}(0)}(x)\,d\tilde{l}\,d\mu_{l_0}(y)\,d\tilde{s}\\
= &\ \frac{C_0e^{C_1\Lambda+C_2K_1T+C_3\sqrt{K_1}r_0}}{r_0^2Vol_{l_0}(B_{l_0}(y_0,r_0))}\frac{\tilde{C}_0e^{\tilde{C}_1\Lambda+\tilde{C}_2K_1T+\tilde{C}_3\sqrt{K_1}r_1}}{r_1^2Vol_{t_0}(B_{t_0}(x_0,r_1))}\\
\quad & \cdot\int^{t_0}_{t_0-4r_0^2}\int_{B_{l_0}(y_0,2r_0)}\int^{l_0+4r_1^2}_{l_0}\int_{B_{t_0}(x_0,2r_1)} H(y,s;x,l)d\mu_{t_0}(x)\,dl\,d\mu_{l_0}(y)\,ds.
\end{aligned}
\end{equation*}

Now let $r_0=r_1=\sqrt{\frac{t_0-l_0}{8}}$, we have by Theorem \ref{davies} that
\begin{align*}
& \ H(y_0,t_0;x_0,l_0)\\
\leq &\ \frac{C_4e^{C_5\Lambda+C_6K_1T+C_7\sqrt{K_1T}}}{(t_0-l_0)^2Vol_{l_0}(B_{l_0}(y_0,\sqrt{\frac{t_0-l_0}{8}}))Vol_{t_0}(B_{t_0}(x_0,\sqrt{\frac{t_0-l_0}{8}}))}\\
\quad & \cdot\int^{t_0}_{\frac{t_0+l_0}{2}}\int^{\frac{t_0+l_0}{2}}_{l_0}\int_{B_{l_0}(y_0,\sqrt{\frac{t_0-l_0}{2}})}\int_{B_{t_0}(x_0,\sqrt{\frac{t_0-l_0}{2}})} H(y,s;x,l)d\mu_{t_0}(x)\,d\mu_{l_0}(y)\,dl\,ds\\
\leq &\ \frac{C_8e^{C_9\Lambda+C_{10}K_1T+C_{11}\sqrt{K_1T}}\sqrt{Vol_{t_0}(B_{t_0}(x_0, \sqrt{\frac{t_0-l_0}{2}}))}\sqrt{Vol_{l_0}(B_{l_0}(y_0, \sqrt{\frac{t_0-l_0}{2}}))}}{Vol_{l_0}(B_{l_0}(y_0,\sqrt{\frac{t_0-l_0}{8}}))Vol_{t_0}(B_{t_0}(x_0,\sqrt{\frac{t_0-l_0}{8}}))}\\
\quad & \cdot exp\left(-\frac{d^2_{t_0}(B_{l_0}(y_0,\sqrt{\frac{t_0-l_0}{2}}),B_{t_0}(x_0,\sqrt{\frac{t_0-l_0}{2}}))}{4e^{4K_1T}(t_0-l_0)}\right)\\
\end{align*}
\begin{align*}
\leq &\ \frac{C_{12}e^{C_{13}\Lambda+C_{14}K_1T+C_{15}\sqrt{K_1T}}}{\sqrt{Vol_{l_0}(B_{l_0}(y_0,\sqrt{\frac{t_0-l_0}{8}}))}\sqrt{Vol_{t_0}(B_{t_0}(x_0,\sqrt{\frac{t_0-l_0}{8}}))}}\\
\quad & \cdot exp\left(-\frac{d^2_{t_0}(B_{t_0}(y_0,e^{K_1T}\sqrt{\frac{t_0-l_0}{2}}),B_{t_0}(x_0,e^{K_1T}\sqrt{\frac{t_0-l_0}{2}}))}{4e^{4K_1T}(t_0-l_0)}\right).
\end{align*}
In the last step above, we have used Lemma \ref{volumecomparison} to get
\begin{align*}
Vol_{t}(B_{t}(z, \sqrt{\frac{t_0-l_0}{2}}))
&\leq\ C_{16}e^{C_{17}\sqrt{K_1T}}\,Vol_{t}(B_{t}(z, \sqrt{\frac{t_0-l_0}{8}}))
\end{align*}
for all $t\in[0,T)$.

Since
\begin{align*}
&d_{t_0}(B_{t_0}(y_0,e^{K_1T}\sqrt{\frac{t_0-l_0}{2}}),B_{t_0}(x_0,e^{K_1T}\sqrt{\frac{t_0-l_0}{2}}))\\
=&\left\{\begin{aligned}
&0,\hspace{5cm}if\  d_{t_0}(x_0,y_0)\leq e^{K_1T}\sqrt{2(t_0-l_0)}\\
&d_{t_0}(x_0,y_0)-e^{K_1T}\sqrt{2(t_0-l_0)},\ \ if\  d_{t_0}(x_0,y_0)>e^{K_1T}\sqrt{2(t_0-l_0)}\end{aligned}\right.,
\end{align*}

it follows that when $d_{t_0}(x_0,y_0)\leq e^{K_1T}\sqrt{2(t_0-l_0)}$,
\begin{align*}
&exp\left(-\frac{d^2_{t_0}(B_{t_0}(y_0,e^{K_1T}\sqrt{\frac{t_0-l_0}{2}}),B_{t_0}(x_0,e^{K_1T}\sqrt{\frac{t_0-l_0}{2}}))}{4e^{4K_1T}(t_0-l_0)}\right)\\
=& 1\\
\leq & \ e^{\frac{1}{4e^{2K_1T}}}exp\left(-\frac{d_{t_0}^2(x_0,y_0)}{8e^{4K_1T}(t_0-l_0)}\right).
\end{align*}

When $d_{t_0}(x_0,y_0)> e^{K_1T}\sqrt{2(t_0-l_0)}$, we have

\begin{align*}
&  exp\left(-\frac{d^2_{t_0}(B_{t_0}(y_0,e^{K_1T}\sqrt{\frac{t_0-l_0}{2}}),B_{t_0}(x_0,e^{K_1T}\sqrt{\frac{t_0-l_0}{2}}))}{4e^{4K_1T}(t_0-l_0)}\right)\\
= &\  exp\left(-\frac{(d_{t_0}(x_0,y_0)-e^{K_1T}\sqrt{2(t_0-l_0)})^2}{4e^{4K_1T}(t_0-l_0)}\right)\\
\leq &\  exp\left(\frac{-1/2\cdot d_{t_0}^2(x_0,y_0)+2e^{2K_1T}(t_0-l_0)}{4e^{4K_1T}(t_0-l_0)}\right)\\
= &\  e^{\frac{1}{2e^{2K_1T}}}\cdot exp\left(-\frac{d_{t_0}^2(x_0,y_0)}{8e^{4K_1T}(t_0-l_0)}\right).
\end{align*}

Therefore, we get

\begin{align*}
& \ H(y_0,t_0;x_0,l_0)\leq \frac{C_{18}exp\left(\displaystyle C_{19}\Lambda+C_{20} KT+C_{21} \sqrt{K_1T}\right)exp\left(\displaystyle -\frac{d_{t_0}^2(x_0,y_0)}{8e^{4K_1T}(t_0-l_0)}\right)}{\sqrt{Vol_{l_0}(B_{l_0}(y_0,\sqrt{\frac{t_0-l_0}{8}}))}\sqrt{Vol_{t_0}(B_{t_0}(x_0,\sqrt{\frac{t_0-l_0}{8}}))}}.
\end{align*}

\end{proof}

The rest of the proof of Theorem \ref{upperbound} follows \cite{CTY2011}. We include them here for the purpose of completeness. The following Lemma shows an $L^1$ bound of the Dirichlet heat kernel.

\begin{lemma}\label{l1lemma}
Let $\Omega$ be a compact manifold with nonempty boundary $\partial \Omega$, and $g(t)$, $t\in[0, T)$, a solution to the Ricci flow \eqref{rf} on $\Omega$. Denote by $H_{\Omega}(y,t;x,l)$ and $G_{\Omega}(x,l;y,t)$ the Dirichlet heat kernel for $\frac{\partial}{\partial t}-\Delta_{t,y}$ and $\frac{\partial}{\partial l}+\Delta_{l,x}-R(x,l)$ on $\Omega$, respectively. Then we have
\begin{equation}\label{l1bound}
e^{-\int_l^t\sup_M R(t)dt}\leq \int_{\Omega} H_{\Omega}(y,t;x,l)d\mu_{t}(y)\leq e^{-\int_l^t\inf_MR(t)dt},
\end{equation}
and
\begin{equation}
\int_{\Omega}G_{\Omega}(x,l;y,t)d\mu_{l}(x)\equiv1
\end{equation}
for any $x, y\in int(\Omega)$.

\end{lemma}

\begin{proof} Since
\begin{align*}
\frac{d}{dt}\int_{\Omega} H_{\Omega}(y,t;x,l)d\mu_{t}(y)&=\int_{\Omega}(\Delta_{t,y}H_{\Omega}(y,t;x,l)-R\cdot H_{\Omega}(y,t;x,l))d\mu_{t}(y)\\
&= \int_{\partial\Omega}\nu_z(H_{\Omega}(z,t;x,l))dS - \int_{\Omega} R\cdot H_{\Omega}(y,t;x,l)d\mu_{t}(y)\\
&= - \int_{\Omega} R\cdot H_{\Omega}(y,t;x,l)d\mu_{t}(y)
\end{align*}
and
\begin{align*}
\frac{d}{dl}\int_{\Omega} G_{\Omega}(x,l;y,t)d\mu_{l}(x)&=\int_{\Omega}-\Delta_{l,x}G_{\Omega}(x,l;y,t)d\mu_{l}(x)\\
&= -\int_{\partial\Omega}\nu_z(G_{\Omega}(z,l;y,t))dS\\
&= 0,
\end{align*}
the Lemma follows immediately.
\end{proof}

Since the heat kernel on a complete manifold is the limit of the Dirichlet heat kernels on a family of exhausting open subsets of the manifold, Lemma \ref{l1lemma} implies that

\begin{corollary}\label{global L1 bound}
Let $(M^n, g_{t})$, $t\in[0,T)$, be a complete solution to the Ricci flow \eqref{rf}. Then we have the following estimates for the heat kernel $H(y,t;x,l)$ and the adjoint heat kernel $G(x,l;y,t)$:
\begin{equation} \label{L1HK}
e^{-\int_l^t\sup_M R(t)dt}\leq \int_{M} H(y,t;x,l)d\mu_{t}(y)\leq e^{-\int_l^t\inf_M R(t)dt},
\end{equation}
and
\begin{equation} \label{L1AHK}
\int_M G(x,l;y,t)d\mu_{l}(x)\equiv 1.
\end{equation}
\end{corollary}

From Corollary \ref{global L1 bound} and the mean value inequality, we can also show the following rough $C^0$ bound of $H$.

\begin{lemma}\label{upperbound2}
Let $(M^n, g(t))$, $t\in[0, T)$ and $T<\infty$, be a complete solution to the Ricci flow. Assume that $Rc\geq -K_1$ for all time $t$, and that $\Lambda=\int_0^T \sup_M|Rc|(t)dt<\infty$. Then there exist constants $\tilde{C}_1,\ \tilde{C}_2,$, $\tilde{C_3}$ and $\tilde{C_4}$ such that
$$H(y,t;x,l)\leq \min\left\{\frac{\tilde{C}_1e^{\tilde{C}_2\Lambda+\tilde{C}_3K_1T+\tilde{C}_4\sqrt{K_1T}}}{\vol_{l}(B_{l}(y,\sqrt{\frac{t-l}{8}}))}, \frac{\tilde{C}_1e^{\tilde{C}_2\Lambda+\tilde{C}_3K_1T+\tilde{C}_4\sqrt{K_1T}}}{\vol_{t}(B_{t}(x,\sqrt{\frac{t-l}{8}}))}\right\}.$$

\end{lemma}

\begin{proof}
In Theorem \ref{mean}, by choosing the parabolic cylinder $P_{g(l)}(y,t,r_0,-(r_0)^2)$ with $r_0=\sqrt{\frac{t-l}{8}}$, we have
\begin{align*}
H(y,t;x,l)&\ \leq \sup_{P_{g(l)}(y,t,r_0,-(r_0)^2)}H(\cdot,\cdot\,;x,l)\\
&\leq\ \frac{C_1e^{C_2\Lambda+C_3K_1T+C_4\sqrt{K_1T}}}{(t-l)Vol_{l}(B_{l}(y,r_0))}\int_{\frac{t+l}{2}}^{t}\int_{B_{l}(y, \sqrt{\frac{t-l}{2}})}H(z,s;x,l)d\mu_{l}(z)ds.
\end{align*}

By \eqref{L1HK}, we can see that
\begin{align*}
\int_{\frac{t+l}{2}}^{t}\int_{B_{l}(y, \sqrt{\frac{t-l}{2}})}H(z,s;x,l)d\mu_l(z)ds&\leq e^{n\Lambda}\int_{\frac{t+l}{2}}^{t}\int_M H(z,s;x,l)d\mu_{s}(z)ds\\
& \leq e^{n\Lambda}\int_{\frac{t+l}{2}}^{t}e^{C_5K_1(s-l)}ds\\
&= e^{n\Lambda}\cdot\frac{e^{C_5K_1(t-l)}-e^{\frac{C_5K_1(t-l)}{2}}}{C_5K_1}.
\end{align*}
Hence,
\begin{align*}
H(y,t;x,l)& \leq \frac{C_1e^{C_2\Lambda+C_3K_1T+C_4\sqrt{K_1T}}}{\vol_{l}(B_{l}(y,\sqrt{\frac{t-l}{8}}))}\cdot \frac{e^{C_5K_1(t-l)}-e^{\frac{C_5K_1(t-l)}{2}}}{C_5K_1(t-l)}\\
&\leq \frac{C_6e^{C_7\Lambda+C_8K_1T+C_9\sqrt{K_1T}}}{\vol_{l}(B_{l}(y,\sqrt{\frac{t-l}{8}}))}.
\end{align*}

Similarly, using \eqref{L1AHK}, we can get
$$H(y,t;x,l)\leq \frac{\tilde{C}_1e^{\tilde{C}_2\Lambda+\tilde{C}_3K_1T+\tilde{C}_4\sqrt{K_1T}}}{\vol_{t}(B_{t}(x,\sqrt{\frac{t-l}{8}}))}.$$
\end{proof}


Now we are ready to prove Theorem \ref{upperbound}.\\

\begin{proof}[Proof of Theorem \ref{upperbound}:]
Let $\sigma=\sqrt{\frac{t_0-l_0}{8}}$ and $r_0=d_{t_0}(x_0,y_0)$.
By Theorem \ref{upperbound1}, we have
\begin{align*}
& \ H(y_0,t_0;x_0,l_0)\leq \frac{C_{0}e^{C_1\Lambda+C_2 K_1T+C_3 \sqrt{K_1T}}exp\left(\displaystyle -\frac{d_{t_0}^2(x_0,y_0)}{8e^{4K_1T}(t_0-l_0)}\right)}{\sqrt{Vol_{l_0}(B_{l_0}(y_0,\sqrt{\frac{t_0-l_0}{8}}))}\sqrt{Vol_{t_0}(B_{t_0}(x_0,\sqrt{\frac{t_0-l_0}{8}}))}}.
\end{align*}

Since $B_{t_0}(x_0,\sigma)\subset B_{t_0}(y_0,\sigma+r_0)\subset B_{l_0}(y_0, e^{\Lambda}(\sigma+r_0))$,
\begin{align*}
\vol^{-1}_{l_0}(B_{l_0}(y_0,\sigma))&\leq \vol^{-1}_{t_0}(B_{t_0}(x_0,\sigma))\cdot \frac{\vol_{t_0}(B_{l_0}(y_0,e^{\Lambda}(\sigma+r_0)))}{\vol_{l_0}(B_{l_0}(y_0,\sigma))}\\
&\leq e^{n\Lambda}\vol^{-1}_{t_0}(B_{t_0}(x_0,\sigma))\cdot \frac{\vol_{l_0}(B_{l_0}(y_0,e^{\Lambda}(\sigma+r_0)))}{\vol_{l_0}(B_{l_0}(y_0,\sigma))}\\
&\leq  e^{n\Lambda}\vol^{-1}_{t_0}(B_{t_0}(x_0,\sigma))\cdot\frac{e^{n\Lambda}(\sigma+r_0)^ne^{\sqrt{(n-1)K_1}e^{\Lambda}(\sigma+r_0)}}{\sigma^n}\\
&\leq e^{n\Lambda}\vol^{-1}_{t_0}(B_{t_0}(x_0,\sigma))\cdot e^{n\Lambda}(1+\frac{r_0}{\sigma})^ne^{C_4\sqrt{K_1T}e^{\Lambda}(1+\frac{r_0}{\sigma})}.
\end{align*}

Let $\hat{C}_1=\frac{1}{64e^{4K_1T}}$, $\hat{C}_2=C_4\sqrt{K_1T}e^{\Lambda}$ and $\eta=\frac{r_0}{\sigma}$, we have
\begin{align*}
& exp\left(\displaystyle -\frac{d_{t_0}^2(x_0,y_0)}{8e^{4K_1T}(t_0-l_0)}\right)\cdot (1+\frac{r_0}{\sigma})^ne^{C_4\sqrt{K_1T}e^{\Lambda}\frac{r_0}{\sigma}}\\
= & exp\left(-\hat{C}_1\eta^2+\hat{C}_2\eta\right)(1+\eta)^n\leq C_5e^{C_6e^{C_7\Lambda+C_8K_1T}}.
\end{align*}
Thus,
\begin{align*}
& \ H(y_0,t_0;x_0,l_0)\leq \frac{C_{9}e^{C_{10}e^{C_{11}\Lambda+C_{12}K_1T}}exp\left(\displaystyle -\frac{d_{t_0}^2(x_0,y_0)}{8e^{4K_1T}(t_0-l_0)}\right)}{Vol_{t_0}(B_{t_0}(x_0,\sqrt{\frac{t_0-l_0}{8}}))}.
\end{align*}


Similarly, one can show that
$$H(y,t_0;x,l_0)\leq \frac{C_{13}e^{C_{14}e^{C_{15}\Lambda+C_{16}K_1T}}exp\left(\displaystyle -\frac{d_{t_0}^2(x,y)}{8e^{4K_1T}(t_0-l_0)}\right)}{Vol_{l_0}(B_{l_0}(y,\sqrt{\frac{t_0-l_0}{8}}))}.$$
\end{proof}

\begin{remark}
If one assumes $Rc\geq0$ on $[0, T)$, and $\Lambda=\int_0^T\sup_M|Rc|(t)dt<\infty$, then the upper bound above can be improved to
$$H(y,t;x,l)\leq C_1e^{C_2\Lambda}exp\left(-\frac{d^2_{t}(x,y)}{8(t-l)}\right)\min\left\{\frac{1}{\vol_{l}(B_{l}(y,\sqrt{\frac{t-l}{8}}))},\frac{1}{\vol_{t}(B_{t}(x,\sqrt{\frac{t-l}{8}}))}\right\}.$$
\end{remark}

\section{A Gaussian lower bound of $H(y,t;x,l)$}

In this section, we prove Theorem \ref{lowerbound} following \cite{CaZh2011}.

\begin{proof}[Proof of Theorem \ref{lowerbound}:]
Let $W$ be a large constant to be determined later. By Theorem \ref{upperbound} and \eqref{L1HK}, we have
\begin{align*}
& \int_{B_{t_0}(x_0,\sqrt{W(t_0-l_0)})}H^2(y, t_0;x_0,l_0) d\mu_{t_0}(y)\\
&\geq \frac{1}{\vol_{t_0}\left(B_{t_0}(x_0,\sqrt{W(t_0-l_0)})\right)}\left(\int_{B_{t_0}(x_0,\sqrt{W(t_0-l_0)})}H(y,t_0;x_0,l_0)d\mu_{t_0}(y)\right)^2\\
&= \frac{1}{\vol_{t_0}\left(B_{t_0}(x_0,\sqrt{W(t_0-l_0)})\right)}\left(\int_M H(y,t_0;x_0,l_0)d\mu_{t_0}(y)\right.\\
&\hspace{6.5cm} \left.-\int_{M-B_{t_0}(x_0,\sqrt{W(t_0-l_0)})}H(y,t_0;x_0,l_0)d\mu_{t_0}(y)\right)^2\\
&\geq \frac{1}{\vol_{t_0}\left(B_{t_0}(x_0,\sqrt{W(t_0-l_0)})\right)}\\
&\quad\  \cdot\left(e^{-C_0\Lambda}-\int\limits_{M-B_{t_0}(x_0,\sqrt{W(t_0-l_0)})}\frac{\tilde{C}_0\hat{C}_0 exp\left(\displaystyle -\frac{d_{t_0}^2(x_0,y)}{8e^{4K_1T}(t_0-l_0)}\right)}{\vol_{t_0}(B_{t_0}(x_0,\sqrt{\frac{t_0-l_0}{8}}))}d\mu_{t_0}(y)\right)^2,
\end{align*}
where $\hat{C}_0=e^{\tilde{C}_1e^{\tilde{C}_2\Lambda+\tilde{C}_3K_1T}}$.

Since
\begin{align*}
& \int_{M-B_{t_0}(x_0,\sqrt{W(t_0-l_0)})}\frac{exp\left(\displaystyle -\frac{d_{t_0}^2(x_0,y)}{8e^{4K_1T}(t_0-l_0)}\right)}{\vol_{t_0}(B_{t_0}(x_0,\sqrt{\frac{t_0-l_0}{8}}))}d\mu_{t_0}(y)\\
\leq &\, e^{-\frac{W}{16e^{4K_1T}}}\int_{M-B_{t_0}(x_0,\sqrt{W(t_0-l_0)})}\frac{exp\left(\displaystyle -\frac{d_{t_0}^2(x_0,y)}{16e^{4K_1T}(t_0-l_0)}\right)}{\vol_{t_0}(B_{t_0}(x_0,\sqrt{\frac{t_0-l_0}{8}}))}d\mu_{t_0}(y)\\
\leq &\, \frac{e^{-\frac{W}{16e^{4K_1T}}}}{\vol_{t_0}(B_{t_0}(x_0,\sqrt{\frac{t_0-l_0}{8}}))}\int_{\sqrt{W(t_0-l_0)}}^{\infty}e^{-\frac{\rho^2}{16e^{4K_1T}(t_0-l_0)}}\frac{d\vol_{t_0}(B_{t_0}(x_0,\rho))}{d\rho}d\rho\\
\end{align*}
\begin{align*}
= &\,\frac{e^{-\frac{W}{16e^{4K_1T}}}}{\vol_{t_0}(B_{t_0}(x_0,\sqrt{\frac{t_0-l_0}{8}}))}\left(-\int_{\sqrt{W(t_0-l_0)}}^{\infty}\vol_{t_0}(B_{t_0}(x_0,\rho))\cdot\frac{d( e^{-\frac{\rho^2}{16e^{4K_1T}(t_0-l_0)}})}{d\rho}d\rho\right.\\
&\hspace{4.5cm}\left.+\,\vol_{t_0}(B_{t_0}(x_0,\rho))e^{-\frac{\rho^2}{16e^{4K_1T}(t_0-l_0)}}\bigg{|}_{\sqrt{W(t_0-l_0)}}^{\infty}\right)\\
\leq&\,
e^{-\frac{W}{16e^{4K_1T}}}\int_{\sqrt{W(t_0-l_0)}}^{\infty}\frac{\vol_{t_0}(B_{t_0}(x_0,\rho))}{\vol_{t_0}(B_{t_0}(x_0,\sqrt{\frac{t_0-l_0}{8}}))}e^{-\frac{\rho^2}{16e^{4K_1T}(t_0-l_0)}}\frac{\rho}{8e^{4K_1T}(t_0-l_0)}d\rho\\
\leq&\, e^{-\frac{W}{16e^{4K_1T}}}\int_{\sqrt{W(t_0-l_0)}}^{\infty}(\frac{\rho}{\sqrt{t_0-l_0}})^ne^{(n-1)\sqrt{K_1}\rho}e^{-\frac{\rho^2}{16e^{4K_1T}(t_0-l_0)}}\frac{\rho}{8e^{4K_1T}(t_0-l_0)}d\rho.
\end{align*}

Set $\eta=\frac{\rho}{\sqrt{t_0-l_0}}$. If we choose $W=C_1e^{C_2\Lambda+C_3K_1T}$ big enough so that
$$\int_{\sqrt{W}}^{\infty}\eta^ne^{(n-1)\sqrt{K_1T}\eta}e^{-\frac{\eta^2}{16e^{4K_1T}}}\frac{\eta}{8e^{4K_1T}}d\eta\leq\frac{1}{2\tilde{C}_0\hat{C}_0} e^{\frac{W}{32e^{4K_1T}}},$$
and
$$e^{-\frac{W}{32e^{4K_1T}}}\leq e^{-C_0\Lambda},$$
then
\begin{align*}
& \int_{M-B_{t_0}(x_0,\sqrt{W(t_0-l_0)})}\frac{\tilde{C}_0\hat{C}_0exp\left(\displaystyle -\frac{d_{t_0}^2(x_0,y)}{8e^{4K_1T}(t_0-l_0)}\right)}{\vol_{t_0}(B_{t_0}(x_0,\sqrt{\frac{t_0-l_0}{8}}))}d\mu_{t_0}(y)\\
\leq & \tilde{C}_0\hat{C}_0e^{-\frac{W}{16e^{4K_1T}}}\int_{\sqrt{W}}^{\infty}\eta^ne^{(n-1)\sqrt{K_1T}\eta}e^{-\frac{\eta^2}{32e^{4K_1T}}}\frac{\eta}{8e^{4K_1T}}d\eta\\
\leq & \frac{1}{2}e^{-\frac{W}{32e^{4K_1T}}}.
\end{align*}

Thus, we can get
\begin{align*}
& \int_{B_{t_0}(x_0,\sqrt{W(t_0-l_0)})}H^2(y, t_0;x_0,l_0) d\mu_{t_0}(y)\\
\geq & \frac{1}{\vol_{t_0}\left(B_{t_0}(x_0,\sqrt{W(t_0-l_0)})\right)}\left(e^{-C_0\Lambda}-\frac{1}{2}e^{-\frac{W}{32e^{4K_1T}}}\right)^2\\
\geq& \frac{1}{4\vol_{t_0}\left(B_{t_0}(x_0,\sqrt{W(t_0-l_0)})\right)}e^{-2C_0\Lambda}.
\end{align*}

By Corollary \ref{global L1 bound}, there exists a point $y_1\in B_{t_0}(x_0, \sqrt{W(t_0-l_0)})$ such that
$$H(y_1,t_0;x_0,l_0)\geq \frac{1}{4\vol_{t_0}\left(B_{t_0}(x_0,\sqrt{W(t_0-l_0)})\right)}e^{-3C_0\Lambda}.$$
From Theorem 3.3 in \cite{Zha2006}, we know
\begin{equation}\label{harnack}
H(y_1,t_0;x_0,l_0)\leq C_3H^{\frac{1}{1+\delta}}(y_0,t_0;x_0,l_0)K^{\frac{\delta}{1+\delta}}e^{2\frac{d^2_{t_0}(y_1,y_0)}{\delta(t_0-l_0)}},
\end{equation}
where $\delta$ is any positive number, and $K=\max_{M\times[\frac{t_0+l_0}{2},t_0]}\{H(y,t;x_0,l_0)\}$.

Since for $t\in[\frac{t_0+l_0}{2},t_0]$, we have
\begin{align*}
\vol_{t}(B_{t}(x_0,\sqrt{\frac{t-l_0}{8}}))&\geq \vol_t(B_t(x_0,\sqrt{\frac{t_0-l_0}{16}}))\\
&\geq e^{-n\Lambda}\vol_{t_0}(B_{t_0}(x_0, e^{-\Lambda}\sqrt{\frac{t_0-l_0}{16}})),
\end{align*}
according to Lemma \ref{upperbound2}, we have for any $(y,t)\in M\times[\frac{t_0+l_0}{2},t_0]$,
\begin{align*}
H(y,t;x_0,l_0)\leq \frac{\tilde{C}_0e^{\tilde{C}_1\Lambda+\tilde{C}_2K_1T+\tilde{C}_3\sqrt{K_1T}}}{\vol_{t}(B_{t}(x_0,\sqrt{\frac{t-l_0}{8}}))}\leq \frac{\tilde{C}_0e^{C_4\Lambda+\tilde{C}_2 K_1T+\tilde{C}_3\sqrt{K_1T}}}{\vol_{t_0}(B_{t_0}(x_0, e^{-\Lambda}\sqrt{\frac{t_0-l_0}{16}}))},
\end{align*}
i.e.,
$$K\leq \frac{\tilde{C}_0e^{C_4\Lambda+\tilde{C}_2 K_1T+\tilde{C}_3\sqrt{K_1T}}}{\vol_{t_0}(B_{t_0}(x_0, e^{-\Lambda}\sqrt{\frac{t_0-l_0}{16}}))}.$$

In \eqref{harnack}, letting $\delta=1$ and noticing that $$d^2_{t_0}(y_1,y_0)\leq 2(d^2_{t_0}(y_1,x_0)+d^2_{t_0}(x_0,y_0))\leq 2W(t_0-l_0)+2d^2_{t_0}(x_0,y_0),$$ we get
\begin{align*}
& H(y_0,t_0;x_0,l_0)\geq\frac{C_5e^{-C_6\Lambda-C_7K_1T-C_8\sqrt{K_1T}}e^{-2W}\vol_{t_0}(B_{t_0}(x_0, e^{-\Lambda}\sqrt{\frac{t_0-l_0}{16}}))}{\vol^2_{t_0}(B_{t_0}(x_0, \sqrt{W(t_0-l_0)}))}e^{-\frac{4d^2_{t_0}(x_0,y_0)}{t_0-l_0}}.
\end{align*}
By Lemma \ref{volumecomparison}, we have
\begin{align*}
\vol_{t_0}(B_{t_0}(x_0, \sqrt{W(t_0-l_0)})
&\leq C_9e^{n\Lambda}W^{\frac{n}{2}}e^{C_{10}\sqrt{WK_1T}}\vol_{t_0}(B_{t_0}(x_0,e^{-\Lambda}\sqrt{\frac{t_0-l_0}{16}})),
\end{align*}
and
$$\vol_{t_0}(B_{t_0}(x_0, \sqrt{W(t_0-l_0)}))\leq C_{11}W^{\frac{n}{2}}e^{C_{12}\sqrt{WK_1T}}\vol_{t_0}(B_{t_0}(x_0,\sqrt{\frac{t_0-l_0}{8}})).$$

Therefore,
\begin{align*}
H(y_0,t_0;x_0,l_0)&\geq\frac{C_{11}e^{-C_{14}e^{C_{15}\Lambda+C_{16}K_1T}}}{\vol_{t_0}(B_{t_0}(x_0, \sqrt{\frac{t_0-l_0}{8}}))}e^{-\frac{4d^2_{t_0}(x_0,y_0)}{(t_0-l_0)}}.
\end{align*}


\end{proof}

\begin{remark}
If $Rc\geq0$ on $[0, T)$ and $\Lambda=\int_0^T\sup_M|Rc|(t)dt<\infty$, then we have
$$H(y,t;x,l)\geq\frac{ C_1e^{-C_2\Lambda}exp\left(-\frac{4d^2_{t}(x,y)}{(t-l)}\right)}{\vol_{t}(B_{t}(x,\sqrt{\frac{t-l}{8}}))}.$$
\end{remark}





\bibliography{mybib}
\bibliographystyle{plain}

\end{document}